\documentclass[MSNbibl,nameyear,dvips]{arxstspdf}
\makeatletter
   \@ifpackageloaded{graphicx}{}{\usepackage{graphicx}}
\makeatother
\usepackage{flushend}
\usepackage{stfloats}
\usepackage{mathbh}
\usepackage{url,breakurl}

\volume{30}
\issue{3}
\pubyear{2015}
\firstpage{391}
\lastpage{412}
\doi{10.1214/15-STS512}% Updated by VTEXPTS2LaTeX.exe, 03.02.2015 11:12
\docsubty{FLA}

\makeatletter
\newcommand{\ind}{\mathbh{1}}
\newcommand{\citett}[1]{\citeauthor{#1}, \citeyear{#1}}
\makeatother

\begin{document}
\begin{frontmatter}
\title{The Ghosts of the \'Ecole Normale}%\thanksref{T1}
% kai straipsnis turi susijusiu diskusiju ir rejoinder'iu
%\relateddois{T1}{Discussed in \relateddoi{d}{10.1214/00-STSXXX} ...;
%rejoinder at \relateddoi{r}{10.1214/00-STSXXXX}.}
\runtitle{The Ghosts of the  \'Ecole  Normale}
%\pdftitle{}
\dedicated{Life, death and legacy of Ren\'e Gateaux\thanksref{T1}}\vspace*{-6pt}

\begin{aug}
% Corresponding author: Laurent MAZLIAK - laurent.mazliak@upmc.fr% Updated by VTEXPTS2LaTeX.exe, 04.02.2015 08:38
%Updated by VTEXPTS2LaTeX.exe, 03.02.2015 11:12
\author[A]{\fnms{Laurent}~\snm{Mazliak}\corref{}\ead
[label=e1]{laurent.mazliak@upmc.fr}}%,
%\author[]{\fnms{}~\snm{}\ead[label=]{}}
% \and
%\author[]{\fnms{}~\snm{}\ead[label=]{}}
\runauthor{L. Mazliak}
%\pdfauthor{}

\affiliation{Universit\'e Pierre et Marie Curie}

\address[A]{Laurent Mazliak is Assistant Professor at Sorbonne Universit\'{e}s, Universit\'{e}
Pierre et Marie Curie, Laboratoire de
Probabilit\'{e}s et Mod\`{e}les Al\'{e}atoires, Case Courrier 188,
4 Place Jussieu, 75252 Paris Cedex 05, France \printead{e1}.}
%\address[]{ \printead{}.}
\end{aug}
\thankstext{T1}{In all the literature, there is a significant
uncertainty regarding whether the name bears a
circumflex accent or not (due to the confusion with the word \textit
{g\^ateau}---cake in French). In the present paper, I shall adopt the
mathematician's own use of \textit{NOT} writing the name with an accent
(this is to conform with his birth certificate).}

% ABSTRACT
%
\begin{abstract}
The present paper deals with the life and some aspects of the
scientific contributions of the mathematician Ren\'e Gateaux, killed
during World War I at the age of 25.
Though he died very young, he left interesting results in functional
analysis. In particular, he was among the first to try to construct an
integral over an infinite-dimensional space. His ideas were extensively
developed later by Paul L\'evy. Among other things, L\'evy interpreted
Gateaux's integral in a probabilistic framework that later contributed
to the construction of the Wiener measure. This article tries to
explain this singular personal and professional destiny in pre- and
postwar France.
\end{abstract}

% KEYWORDS
\begin{keyword}
\kwd{History of mathematics}
\kwd{functional analysis}
\kwd{integration}
\kwd{probability}
\kwd{Wiener measure}
\end{keyword}
\end{frontmatter}

\setcounter{footnote}{1}

%s1 #&#
\section{Introduction}\label{sec1}

In his\vspace*{1pt} seminal 1923 paper on Brownian motion, Norbert Wiener
mentioned\footnote{\citet{Wiener1923}, page 132.} that \textit{integration
in infinitely many dimensions \textup{(}was\textup{)} a relatively little-studied
problem} and that \textit{all that has been done on it \textup{(}was\textup{)} due to
Gateaux\textup{,} L\'evy\textup{,} Daniell} and himself. Following Wiener, the most
complete investigations had been those \textit{begun by Gateaux and
carried out by L\'evy}.

It was in 1922 that L\'evy's book \textit{Le\c cons d'Analyse
Fonctionnelle} (\citett{Levy1922}) was published after his lectures given
at the Coll\`ege de France in the aftermath of the Great War. L\'evy's
book, and, more specifically, L\'evy himself, made a profound
impression on Wiener. The American
mathematician emphasized how L\'evy explained personally to him how his
own method of integration in infinitely many dimensions, which extended
results L\'evy found in\vadjust{\goodbreak} Gateaux's works, was the convenient tool he
needed for his construction of Brownian motion measure.

I shall comment later on the path linking Gateaux's works to L\'evy's
fundamental studies, but let me begin by discussing the circumstances
which constituted the initial motivation behind the current paper.
Gateaux was killed at the very beginning of the Great War in October of
1914. He died at the age of 25, before having obtained any academic
position, even before having completed a doctorate. His publications
formed a rather thin set of a few notes presented to the Academy of
Sciences of Paris and to the Accademia dei Lincei of Rome. None of them
dealt with infinite-dimensional integration. Neverthess, Gateaux's name
is still known today, and even to (some) undergraduate students,
through a basic notion of calculus known as {\it Gateaux
differentiability}.\footnote{Let me recall that Gateaux
differentiability of a function $\phi$ defined on $\mathbb{R}^n$ is the
directional differentiability: $\phi$ is said to be Gateaux
differentiable at $\theta\in\mathbb{R}^n$ if for any vector $h$
given in $\mathbb{R}
^n$, the function $t \mapsto\phi( \theta+t h)$ has a derivative at 0.
Various notions of differentiability for a function have been
considered by several French mathematicians under Volterra and
Hadamard's influence during the first half of the 20th century. In the
1920s, Hadamard introduced an intermediate concept between Fr\'echet
and Gateaux differentiability. In modern terminology, a function
$\phi\dvtx
E\rightarrow F$, where $E$ and $F$ are two normed spaces, is Hadamard
differentiable at $\theta\in E$ if there is a continuous linear
function $\phi_\theta'\dvtx E\rightarrow F$ such that, for any $h\in E$
and any choice of a family $(h_t)_{t >0}$ in $E$ such that
$h_t\rightarrow h$, one has
\[
\lim_{t\rightarrow0} \biggl\|\frac{\phi(\theta+ th_t)-\phi(\theta)}{t}-\phi_\theta'(h)
\biggr\|_F=0.
\]
The difference between Gateaux and Hadamard differentiability is that,
for the latter, the direction $h_t$ is allowed to change in the ratio.
On this topic see \citet{Barbut2014}, Section~4.2, pages 15--17.
Hadamard-differentiabilty is in particular adequate to deal with some
asymptotic estimates in Statistics (see, e.g., \citett{Vaart1998}, Chapter~20---especially page 296 and seq.).} The notion,
weaker than the (now) classical Fr\'echet differentiability, was
mentioned in Gateaux's note (\citett{Gateaux3}, page  311), under the
name \textit{variation premi\`ere} of a functional, though it was probably
already considered by him in 1913 as the name appears in \citet{Gateaux1}, page~326, but without any definition. Regardless, this
notion was in fact only a technicality introduced by Gateaux among the
general properties that a functional can have. L\'evy was probably the
first to name it after Gateaux.\footnote{In \citet{Levy1922}, page 51,
under the name \textit{diff\'erentielle au sens de Gateaux}, Sanger (\citeyear{Sanger1933}) compared the various definitions formulated for the
differential of a functional in his survey about Volterra's functions
of lines. See, in particular, Chapter II on pages 240--253. Gateaux's
definition is considered on pages 250--251.}

So, I wanted to understand how a basic notion of calculus had been
given the name of an unknown mathematician, who died so young, before
having obtained any academic position and even before having defended a
thesis. Such a paradox deserved to be unraveled. It is this apparent
contradiction that I want to address in this paper by presenting Ren\'e
Gateaux's life and death, some of his mathematical research and the
path explaining why we still remember him though so many of his fellows
killed during the war became only a \textit{golden word on our public
squares}, following Aragon's beautiful expression.\footnote{D\'ej\`a la
pierre pense o\`u votre nom s'inscrit\\ D\'ej\`a vous n'\^etes plus
qu'un mot d'or sur nos places\\ D\'ej\`a le souvenir de vos amours
s'efface\\ D\'ej\`a vous n'\^etes plus que pour avoir p\'eri (\citett{Aragon1956}).}

Let me immediately reveal the key to our explanation. Beyond his tragic
fate, Gateaux had two strokes of good fortune. The first one was
related to the main mathematical theme he was interested in, \textit{Functional Analysis} (Analyse Fonctionnelle) in the spirit of Volterra
in Rome and Hadamard in Paris, often also called by them
\textit{functional calculus} (calcul fonctionnel).\footnote{In the sequel, I
shall use the expression \textit{functional analysis} only in reference to
the theories initiated by Volterra, though it today has a slightly
different meaning.} At the beginning of the 20th Century, this subject
was still little studied. In the years following World War I, it
received unexpected developments, in particular, in the unpredictable
direction of probability theory. Gateaux was therefore posthumously in
contact with a powerful stream leading to the emergence of some central
aspects of modern probability, such as Brownian motion as we have seen
in Wiener's own words. It is very fortunate for the historian that
important archival documents about Gateaux's beginnings in mathematics
are still available. Gateaux had in particular been in correspondence
with Volterra before, during and (for some weeks) after a sojourn in
Rome with the Italian mathematician. His letters still exist today at
the \textit{Accademia dei Lincei} and provide precious insight into
Gateaux's first steps. Letters exchanged between Borel and Volterra
about the young man's projects and progress are also available. One
such document is a letter from Gateaux to Volterra dated from 25 August
1914 and written on the battlefield. Moreover, some other material is
accessible such as the military dossier, some of Gateaux's own drafts
of reports about his work, and some scattered letters from him or about
him by other people. This allows us to attempt to reconstruct the life
of the young mathematician during his last seven or eight years.

But it is mainly due to the second stroke of fortune that some memory
of Gateaux (or, at least, of his name) was preserved. Before he went to
the war, Gateaux had left his papers in his mother's house. Among them
were several half-completed manuscripts which were intended to become
chapters of his thesis. After the death of her son, his mother sent the
papers to the \'Ecole Normale. Hadamard collected them and in 1919 passed
them to Paul L\'evy in order to prepare an edition in Gateaux's honor.
Studying Gateaux's papers came at a crucial moment in L\'evy's career.
Not only did they inspire L\'evy's book (\citett{Levy1922}), but they were a
major source for his later achievements in probability theory.

The aim of the present paper is twofold: one aspect is to present an
account of Gateaux's life by using valuable new archival material
discovered in several places, the other is to give some hints of how
his works were completed and---considerably---extended by L\'evy. In
that respect, it is clear that the mathematical ideas of Gateaux were
developed in a direction he could not have expected; probability, for
instance, was absolutely not in his mind. The appearance of the
mathematics of randomness in this inheritance is undoubtedly entirely
due to L\'evy's powerful imagination. It is therefore well beyond the
scope of this article, centered on Gateaux, to present a detailed study
of L\'evy or Wiener's studies on Brownian motion. The interested reader
may refer to several historical expositions such as \citet{Kahane1998}
or \citet{Barbut2014}, pages  54--60. An account from direct
participants in this story can be found in L\'evy's autobiography
(\citett{Levy1970}, page~96 and seq.), or It\^o's comments on Wiener's papers
(\citett{WienerOC}, pages 513--519).

Looking backward, Gateaux's role must not be overestimated in the
history of mathematics. Contrary to some other examples of
mathematicians who died young, such as Abel to cite a famous example,
Gateaux had not made decisive progress in any important direction. So
maybe some words are necessary to explain what a biographical approach
of someone like Gateaux can teach us. The main point here is related to
the Great War and the effect it produced on French mathematicians.

In her memoirs (\citett{Marbo1967}), written at the end of the 1960s, the
novelist Camille Marbo,\footnote{Marbo is Marguerite Appel's nom de
plume. She was the daughter of the mathematician Paul Appell.} Emile
Borel's widow, mentioned that after the end of World War I, her husband
declared that he could not bear any more the atmosphere of the \'Ecole
Normale in mourning, and decided to resign from his position of Deputy
Director. In 1910 Borel had succeeded Jules Tannery in the position,
during a time of extraordinary success for Analysis in France with
outstanding mathematicians such as Henri Poincar\'e, Emile Picard,
Jacques Hadamard, Henri Lebesgue and naturally Borel himself.

A superficial, though impressive, picture of the effect of WWI on the
French mathematical community is read through the personal life of the
aforementioned mathematicians---with the obvious exception of Poincar\'e who had died in 1912. Picard lost
one son in 1915, Hadamard two sons
in 1916 (one in May, one in July) and Borel his adopted son in 1915.
The figures concerning casualties among the students of the \'Ecole
Normale, and especially among those who had just finished their three
year studies at the \textit{rue d'Ulm}, are terrible.\footnote{They were
collected in a small brochure published by the \'Ecole Normale at the end
of the war (\citett{ENSguerre}).} Out of about 280 pupils who entered the
\'Ecole Normale in the years 1911 to 1914, 241 were sent to the front
directly from the school and 101 died during the war. If the President
of the Republic Raymond Poincar\'e could declare that \textit{the \'Ecole of
1914 has avenged the \'Ecole of 1870},\footnote{L'\'{E}cole de 1914 a veng\'e
l'\'Ecole de 1870 (\citett{ENSguerre}, page~3).} the price to pay had
been so enormous that it was difficult to understand how French science
could survive such a hemorrhage. Most of the vanished were brilliant
young men, expected successors of the brightest scholars from the
previous generation in every domain of knowledge. They were so young
that almost none had time to start making a reputation of his own
through professional achievement. As testimony of his assumed
abnegation, Fr\'ed\'eric Gauthier, a young hellenist, who had entered
the \'Ecole Normale Sup\'erieure in 1909 and was killed in July 1916 in
the battle of Verdun, left a melancholic comment on this time of
resignation: \textit{My studies\textup{,} it is true\textup{,} will remain sterile\textup{,} but my
ultimate actions\textup{,} useful for the country\textup{,} have the same value as a
whole life of action}.\footnote{Mes \'etudes, il est vrai, seront
demeur\'ees st\'eriles, mais mes actions derni\`eres, utiles au pays,
vaudront toute une vie d'action (\citett{AnnuaireENS}).}

Gateaux, who died at the very beginning of the war, appears therefore
to be a good representation of the lost generation of \textit{normaliens}
that I have just mentioned; he was at the same time an exception, as
his very name, contrary to almost every one of his companions of
misfortune, was retained in mathematics. The way in which it was
retained and, above all, the direction in which his works received
their most important development [Wiener's seminal paper (\citett{Wiener1923})] was, at least partly, related to the war. L\'evy wrote to
Fr\'echet in 1945:
\begin{quote}
As for myself, I learned the first elements of probability
during the spring of 1919 thanks to Carvallo (the director of studies
at the \'Ecole Polytechnique) who asked me to make three lectures on that
topic to the students there. Besides, in three weeks, I succeeded in
proving new results. And never will I claim for my work in probability
a date before 1919. I can even add, and I told M. Borel so, that I had
not really seen before 1929 how important were the new problems implied
by the theory of denumerable probabilities. But I was prepared by
functional calculus to the studies of functions with an infinite number
of variables and many of my ideas in functional analysis became without
effort ideas which could be applied in probability.\footnote{\citet{Barbut2014}, page 139.}
\end{quote}

The urgent need to renew the teaching of probability at the \'Ecole
Polytechnique was a side effect of the war, when much probabilistic
technique had been used to direct artillery. And it is because Gateaux
was dead that L\'evy was in possession of his papers. Nobody can tell
what L\'evy's career would have been without the conjunction of these
two disparate elements that his fertile mind surprisingly connected.

I began to be interested in Gateaux's story when we were preparing the
edition of Fr\'echet and L\'evy's correspondence with Barbut and Locker
in 2003 [an English edition (\citett{Barbut2014}) was recently published].
Since that date, a lot of work has been done concerning the involvement
of scientists in the Great War, resulting in an increasing number of
publications, and, in particular, the approach of the centennial year
was met by a flow of papers and books in many countries so that it is
difficult to provide an exhaustive list. Let me mention, among many
others, the interesting contributions [\citet{Pepe2011},
\citet{Onghena2011} or the books \citet{AubinGoldstein2014} and
\citet{Downing2014}]. By the way, the centennial was also an occasion for
economists to remember Gateaux's work (\citett{DuggerLambert2013}).

A focus on Gateaux therefore allows us to shed some light on some
specific aspects of mathematics before and after the Great War and to
understand how such an event may have influenced their development, not
only in technical aspects but also because of its terrible human cost.

The paper is divided into four parts. In the first I describe Gateaux's
life before he went to Rome in 1913. Then I present the critical period
in Rome with Volterra. The third part treats his departure to the army
and his last days. Finally, there is a slightly more technical part
which considers the work of Gateaux and how it was recovered by L\'evy
and considerably extended by him so that it became a step toward the
construction of an abstract integral in infinite dimensions and then of
modern probability theory.

%s2 #&#
\section{A Provincial in Paris}

We do not know much about Gateaux's life before he entered the \'Ecole
Normale. Gateaux did not belong to an important family and, moreover,
his family unit consisted only of his parents, his younger brother
Georges and himself. Neither of the brothers had direct descendants, as
both boys died during WWI. I have met a distant member of his family,
namely, the great-great-great-great-grandson of a
great-great-great-grandfather of Ren\'e Gateaux, Mr Pierre Gateaux, who
still lives in Vitry-le-Fran\c{c}ois and most kindly offered access to
the little information he has about his relative.

Ren\'e Eug\`ene Gateaux was born on 5 May 1889 in Vitry-le-Fran\c{c}ois
in the d\'epartement of Marne, 200 km east from Paris.\footnote{Abraham de Moivre was born there 222 years earlier, before the wars of religion
forced him to leave for London where he spent all his scientific
career. Fran\c cois Jacquier was also born there 178 years earlier. A
local historian from Vitry, Gilbert Maheut, has written several short
papers about his three mathematician fellow-citizens. See, in
particular, \citet{Maheut2000}.} Ren\'e's father Henri, born in 1860,
was a small contractor who owned a saddlery and cooperage business in
the outskirts of Vitry. His mother was Marie Roblin, born in Vitry in
1864. Ren\'e's family on his father side came originally from the small
town of Villers-le-Sec at 20 km from Vitry, the rural nest of Gateaux's
family. Ren\'e's birth certificate indicates that Eugene Gateaux
(Henri's father) was a proprietor and Jules Roblin (Marie's father) was
a cooper; the grandparents acted as witnesses when the birth was
registered at the town hall. Eug\`ene's birth certificate indicates
that he was born in Villers-le-Sec in 1821 and that his father was a
carpenter. Perhaps Ren\'e's grandfather came to Vitry to create his
business and employed Marie's father as a cooper. As already mentioned,
the couple had two children: Ren\'e is the elder; the second one,
Georges, was born four years later in 1893. Ren\'e's father died young,
in 1905, aged 44, and the resulting precarious situation may have
increased the boy's determination to succeed in his studies.

I have no details on Ren\'e Gateaux's school career; he was a pupil in
Vitry and then in Reims. The oldest handwritten document I have found
is a letter to the Minister of Public Instruction on 24 February 1906
asking for permission to sit for the examination for admission to the
\'Ecole Normale Sup\'erieure\footnote{After the defeat of 1870, the
prestige of the \'Ecole Polytechnique faded and the \'Ecole Normale Sup\'erieure became the major center of scientific life in France at the
turn of the century. The \'Ecole Polytechnique was to regain a real
importance for scientific research only much later in the 20th Century.
Paul L\'evy, who chose to go to the Polytechnique instead of the \'Ecole
Normale to please his father, was a real exception in mathematical
research at the beginning of 20th Century. He also slightly suffered
from the situation by not belonging to Borel's or Hadamard's usual
network of \textit{normaliens}.} (science division), although he had not
reached the regular minimum age of 18.

Two things can be deduced from this document. The first is that Gateaux
was\vadjust{\goodbreak} a student in a Classe Pr\'eparatoire in the lyc\'ee of
Reims.\footnote{The Classes Pr\'eparatoires are the special sections in
the French school system that train students for the competitive
examinations for entry to the ``Grandes \'Ecoles,'' such as the \'Ecole
Polytechnique or the \'Ecole Normale Sup\'erieure.} Our second
inference is that Gateaux was a brilliant student in his science
classes. He probably obtained his baccalaur\'eat in July 1904 at the
age of only 15. Gateaux was a sufficiently exceptional case for an
inspector (coming at the Lyc\'ee of Reims in March 1907) to mention in
his report that Gateaux had obtained the extraordinary mark of 19 (out
of 20) to a written test in mathematics.\footnote{Archives d\'epartementales de la Marne.} However, he was not admitted to the \'Ecole
Normale on his first attempt in 1906, but only in October 1907 after a
second year in the class of Math\'ematiques Sp\'eciales, as was usually
the case.

What was it like to be a provincial in Paris? Jean Gu\'ehenno, born in
1890, and admitted in 1911 in the literary section, has written some
fine pages on the subject in his {\it Journal d'un homme de 40 ans} (\citett{Guehenno1934}---see, in particular, Chapter VI,
``Intellectuel''). There he describes the \'Ecole Normale Sup\'erieure of
the years before the Great War through the eyes of a young man from a
poor provincial background (much poorer, in fact, than Gateaux's) and
how he was dazzled by the contrast between the intellectual riches of
Paris and the laborious tedium of everyday life in his little
industrial town in Brittany. We also have an obituary (\citett{AnnuaireENS}, pages 136--140) written in 1919 by two of Gateaux's
fellow students from the 1907 science section of the \'Ecole Normale,
Georges Gonthiez and Maurice Janet. They described Gateaux as a good
comrade with benevolence and absolute sincerity, who soon appeared to
his fellow students as one of the best mathematicians of the group.

After the entry at the \'Ecole occurred an event in the young man's life
of undoubted importance since Gonthiez and Janet devote many lines to
it. Gateaux became a member of the Roman Catholic Church. He joined the
church \textit{with fervour}, wrote his two fellows. Such a decision in
1908 may seem surprising: the separation laws between Church and State
had been passed in 1905 and the Roman Church stood accused for its
behavior during the Dreyfus Affair. However, there was concurrently a
revival of interest in Catholicism as a counterweight to triumphant
positivism. Such a current was well represented at the \'Ecole Normale
(\citett{Gugelot1998}). Among Gateaux's fellows was Pierre Poyet, who chose
a religious life and died a few months before he could make his vows as
a Jesuit.

Ren\'e's conversion to Catholicism, which had a profound effect on his
spiritual life, created difficulties for him at the \'Ecole Normale.
Gateaux explained in a letter to Poyet (quoted in \citett{Bessieres1934})
that his conversion was received badly by his fellows and some
professors. Several pages are devoted to Gateaux in B\'essi\`eres'
biography of Poyet (\citett{Bessieres1934}). So far, all efforts to
locate Poyet's personal papers have been fruitless, nevertheless, the
obituary by Gonthiez and Janet in \citet{AnnuaireENS} testifies not
only to the incomprehension felt by Gateaux's fellows, but also to how
they were impressed by the similarity of the methods used by him to
progress in his mathematical and spiritual lives.\footnote{\citet {Bessieres1934} provides a surprising picture of the mystic atmosphere
present at the \'Ecole Normale around Poyet.}

In 1910, Gateaux passed the Agr\'egation of mathematical sciences where
he obtained the 11th rank out of~16. This was not a very good rank, so
it left him no possibility of obtaining a grant to devote himself
entirely to research, as had been the case for Joseph P\'er\`es, for
instance (on which I shall comment later). On 8 July 1912, a
ministerial decree appointed Gateaux as Professor of Mathematics at the
Lyc\'ee of Bar-le-Duc, the principal town of the d\'epartement of
Meuse, 250 km east from Paris, and not very distant from his native town.

Before taking up this position, Gateaux should have fulfilled his
military obligations. From March 1905 (\citett{JO1905}), a new law
replaced the July 1889 regulation for the organization of the army. The
period of active military service had been reduced to 2 years, but
conscription became, in theory, absolutely universal. Gateaux was
particularly affected by article 23 stipulating that the young men who
entered educational institutions such as the \'Ecole Normale Sup\'erieure
could, at their choice, fulfill the first of their two years of
military service in the ranks before their admission to these
institutions or after their exit. Gateaux had chosen the latter option
(\citett{Gateauxarchmil}). In October 1910, Gateaux joined the 94th
Infantry regiment where he was a private. In February 1911 he was
promoted to {\it caporal} (corporal), and finally was declared second
lieutenant in the reserve in September 1911. He had to follow some
special training for officers; the comments made by his superiors on
the military file indicate that the supposed military training at the
\'Ecole Normale had been more virtual than real. On the special pages
devoted to his superior's appraisal, one reads that, though having very
good spirit, Ren\'e was hardly prepared for his rank, but that the
second semester 1912 (which ended in fact in September 1912) seems to
have been better. He had followed a period of instruction for shooting
and obtained very good marks. A final comment in the military file has
a strange resonance with what happened two years later. Gateaux's
superior mentioned that he was able to lead a machine-gun section.

In October 1912, Gateaux, freed from the active army, began his
lectures at the Lyc\'ee of Bar-le-Duc. Gateaux's (very thin) personnel
file contains a personal identification form and a decree of the
Minister of Public Instruction on 2 October 1913 granting him
one-year's leave with an allocation of 100 francs for that year, as
well as a handwritten document showing that he had obtained a David
Weill grant for an amount of 3000 francs.\vspace*{-3pt}

%s3 #&#
\section{The Roman Stage}

Gateaux had indeed begun to work on a thesis with themes closely
related to functional analysis \textit{\`a la Hadamard}. I have found no
precise information about how Gateaux chose this subject for his
research, but it is plausible that he was advised to do so by Hadamard
himself. In 1912, Hadamard had just delivered a series of lectures on
functionals at the Coll\`ege de France and had entered the Academy of
Science in the same year. Paul L\'evy had, moreover, defended his own
brilliant thesis on similar questions in 1911. As well, a young French
normalien of the year before Gateaux, Joseph P\'er\`es, had in
1912--1913 benefited from a David Weill grant offered for a one-year
stay in Rome with Volterra. Volterra himself, invited by Borel and
Hadamard, had come to Paris for a series of lectures on functional
analysis, edited by P\'er\`es and published in 1913 (\citett{Volterra1913}). These were thus good reasons for Gateaux to be
attracted by this new and little explored domain. For a young doctoral
student the natural people to be in contact with were Hadamard in Paris
and Volterra in Rome.\footnote{On Hadamard, a star of the French
mathematical stage of the time, the reader can refer to the book (\citett{MazyaShaposhnikova1998}). Two biographies of Vito Volterra have
recently been published (\citett{Goodstein2007}; \citett{GuerragioPaoloni2008}), and the reader can also find information in
the annotated edition of the correspondence between Volterra and his
French colleagues during WWI (\citett{MazliakTazzioli2009}).} P\'er\`es's example encouraged Gateaux to go to Rome. Some years later, when
Hadamard wrote a report recommending Gateaux for the posthumous
attribution of the Franc\oe ur prize, he mentioned that the young man\vadjust{\goodbreak}
had been \textit{one of those who\textup{,} inaugurating a tradition that could not
be overestimated\textup{,} went to Rome to become familiar with M. Volterra's
methods and theories}.\footnote{Il fut un de ceux qui, inaugurant une
tradition \`a laquelle nous ne saurions trop applaudir, all\`erent \`a
Rome se former aux m\'ethodes et aux th\'eories de M. Volterra (\citett{Hadamard1916}). On the development of student exchanges between Paris
and Rome in these years, see \citet{Mazliak2015}.}

On the occasion of the centennial of Volterra's birth, in 1960, a
volume was edited by the Accademia dei Lincei in Rome in which Giulio
Krall devoted several pages to Volterra's research on the phenomenon of
\textit{hysteresis}, the ``memory of materials,'' which describes the
dependence on time of the state of deformation of certain materials. To
model such a situation, Volterra was led to consider \textit{functions of
lines} (funzione di linea), later called \textit{functionals}
(fonctionnelle) by Hadamard and his followers, which is to say a
function of a real function representing the state of the material, and
to study the equations they must satisfy. These equations happen to be
an infinite-dimensional generalization of partial differential
equations. As Krall mentions,\footnote{\citet{Krall1961}, page~17.}
\textit{from mechanics to electromagnetism\textup{,} the step was small}, and Volterra's
model was applied to different physical situations, such as
electromagnetism or sound produced by vibrating bars.\footnote{Volterra himself was involved in this subject through an important collaboration
with Arthur Gordon Webster from Clark University in the USA. See the
interesting webpage \url{http://physics.clarku.edu/history/history.html\#webster}.} In 1904,
the King made Volterra a Senator of the Kingdom, mostly honorary, but
giving the recipient some influence through his proximity with the men
of power.

Such a combination of science and politics appealed to Borel, who had a
deep friendship with Volterra.\footnote{On the beginning of the
relationship between Borel and Volterra, see \citet{Mazliak2015}.}
Borel had a part in Gateaux's decision to go to Rome, at least as an
intermediary between the young man and Volterra. We indeed find a first
indication of this Roman project in their correspondence. Borel wrote
to Volterra on 18 April 1913 that he intended to support Ren\'e's
request for the grant,
and joined a letter written by Ren\'{e} Gateaux where he explained his research agenda.
Borel then asked Volterra to write a short
letter of support for this project to Liard, the Vice-Rector of the
Paris Academy, and also mentioned that he lent the two books published
by Volterra on the functions of lines to Gateaux [more precisely, the
book \citet{Volterra1912} and the proofs of \citet{Volterra1913}]. On 30
June 1913, Borel communicated the good news to Volterra: a David Weill
grant had been awarded to Gateaux for the year 1913--1914.

Gateaux's aforementioned letter to Borel\footnote{Dated from
Bar-le-Duc, 12 April 1913.} was in Volterra's archives, and,
consequently, we know precisely what his mathematical aims were when he
went to Rome. Gateaux considered two main points of interest for his
future research. The first one is classified as \textit{Fonctionnelles
analytiques} (Analytical functionals) and is devoted to the extension
of the classical results on analytical functions: the Weierstrass
expansion, the equivalence between analyticity and holomorphy and the
Cauchy formula. The second one is devoted to the problem of integration
of a functional.

Gateaux started from the definition Fr\'echet had proposed in 1910 for
an analytical functional (\citett{Frechet1910}) based on a generalization\vspace*{1.5pt}
of a Taylor expansion. A functional\footnote{Throughout the paper, the
functionals considered are always defined on the set of real functions
over a given interval $[a,b]$.} $U$ is homogeneous with order $n$ if
for any $p\ge1$ and any given continuous functions $g_1,\ldots, g_p$
over $[a,b]$, the function defined on $\mathbb{R}^p$ by
\[
(\lambda_1, \dots, \lambda_p)\mapsto U(
\lambda_1g_1+\cdots+ \lambda_pg_p)
\]
is a homogeneous polynomial of degree less than $n$.\footnote{Fr\'echet's definition is in fact given in a different way by means of a
property inspired by a characterization he had proved for real
polynomials (\citett{Frechet1910}, page 204); however, he proves (page~205) that the two properties are equivalent.} Now, a functional $U$ is
by definition analytical if it can be written as
\[
U(f)=\sum_{n=0}^\infty U_n(f),
\]
where $U_n$ are homogeneous functionals of order $n$ (\citett{Frechet1910}, page 214; see also Taylor, \citeyear{Taylor1970}).

Gateaux first proposed to obtain properties of the terms $U_n(f)$ in
the previous expansion of an analytical functional. Then, he intended
to obtain the equivalence between the analyticity of the functional $U$
and its complex differentiability (holomorphy) and to deduce a
definition of analyticity by a Cauchy formula. For that purpose, as he
wrote, one needs a definition of the integral of a real continuous
functional over a real functional field. This may be the first
appearance of questions around infinite-dimensional integration. In
this programmatic letter, Gateaux suggested the way he wanted to
proceed, inspired by Riemann integration:
\begin{quote} Let us restrict ourselves to the definition of the
integral of $U$ in the field of the functions $0\le f\le1$. Let us
divide the interval $(0,1)$ into $n$ intervals. $(\ldots)$ Consider next
the function $f$ in any of the partial intervals as equal to the
numbers $f_1,\ldots,f_n$ which are between 0 and 1. $U(f)$ is a
function of the $n$ variables $f_1\cdots f_n \dvtx  U_n(f_1,\ldots,f_n)$. Let
us consider the expression
\begin{eqnarray*}
I_n&=&\int_0^1\!\int
_0^1\cdots\int_0^1
U_n(f_1,\ldots,f_n)\\
&&\hspace*{16pt}\qquad\qquad{}\cdot df_1\cdots
\,df_n.
\end{eqnarray*}
Suppose that $n$ increases to infinity, each interval converging to 0,
and that $I_n$ tends to a limit $I$ independent of the chosen
divisions. We shall say that $I$ is the integral of $U$ over the field
$0\le f\le1$.
\end{quote}
Gateaux's intention was to study whether the limit $I$ exists for any
continuous functional $U$ or if an extra hypothesis was necessary. In
the last paragraph, Gateaux mentioned the possible applications of this
integration of functionals, such as the residue theorem. All the
applications he mentioned belong in addition to the theory of functions
of a line. There is no hint of a possible connection with potential
theory. This does not appear in the papers published by Gateaux. As it
is a central theme of Gateaux's posthumous texts, it is plausible that
he became conscious of the connection only during his stay in Rome---perhaps under Volterra's influence.

On 28 August 1913, Gateaux wrote directly to Volterra for the first
time, informing him of his arrival in October and also mentioning that
he had already obtained several results for the thesis in Functional
Analysis which he was working on. Gateaux may have enclosed a copy of
his first note to the Comptes-Rendus (\citett{Gateaux1}), published on 4
August 1913 and containing the beginning of his proposed program. The
note is in fact rather limited to an exposition of results and does not
contain any proof, apart from a sketch of how to approximate a
continuous functional $U$ by a sequence of functionals of order $n$
uniformly over each compact subset of the space of continuous real
functions on $[0,1]$.\footnote{We need not dwell upon this technical
result here, which had already been obtained by Fr\'echet previously
(\citett{Frechet1910}, page 197) in a slightly more intricate way. Let me
only observe that Gateaux's elementary technique involves the
replacement of the function $z$ by a linear function over each
subdivision $[\frac{i}{n}, \frac{i+1}{n}]$ of the interval $[0,1]$. A final
perfecting of Gateaux's proof is presented by L\'evy (\citeyear{Levy1922}), pages 105--107.}

About Gateaux's stay in Rome, I do not have many details. An
interesting document, found in the Paris Academy, is the draft of a
report written by Gateaux at the end of his stay for the David Weill
foundation.\footnote{A very touching aspect of the report written by
Gateaux for the David Weill foundation can be found in the pages where
he described the nonmathematical aspects of his journey. Gateaux
mentioned how he regretted that Italy and the Italian language were so
little known in France, when, on the contrary, France and French were
widely known within Italian society.} He mentioned there that he had
arrived in Rome in the last days of October and that he followed two of
Volterra's courses in Rome (one in Mathematical Physics, the other
about application of functional calculus to Mechanics). Gateaux seems
to have worked quite actively in Rome. A first note to the Accademia
dei Lincei (\citett{Gateaux2}) where he extended the results of his
previous note to the Paris Academy was published in December 1913. On a
postcard sent by Borel to Volterra on 1 January 1914, Borel mentioned
how he was glad to learn that Volterra was satisfied with Gateaux. The
young man published three more notes during his stay (\citeauthor{Gateaux3}, \citeyear{Gateaux3}, \citeyear{Gateaux4}, \citeyear{Gateaux5}), and also began to write
more detailed articles---found after the war among his papers.\footnote{L\'evy (in \citett{Gateaux7}, page 70) mentioned that, in one case, two
versions of the same paper were found, both dated March 1914.}

On 14 February 1914, Gateaux made a presentation to Volterra's
seminar\footnote{His lecture notes were found among his papers.} in
which he mainly dealt with the notion of functional differentiation. He
recalled that Volterra introduced this notion to study problems
including hereditary phenomena, and also that it was used by others
(Hadamard and Paul L\'evy) to study some problems of mathematical
physics---such as the equilibrium problem of fitted elastic plates---through the resolution of equations with functional derivatives.

Gateaux came back to France at the beginning of the summer, in June
1914. He expected to go back soon to Rome, as he was almost certain, as
Borel had written to Volterra,\footnote{Borel to Volterra. 3 April
1914.} to obtain the Commercy grant he had applied for. Gateaux soon
wrote that the grant had been awarded.\footnote{Gateaux to Volterra, 14
July 1914.} In the same letter, he mentioned that he had completed a
first version of a note on functionals requested by Volterra to append
it to the German translation of his lectures on functions of lines
(\citett{Volterra1913}). During this month, he had also met the Proviseur
of the Lyc\'ee in Bar-le-Duc on July 20th, as the man sadly observed in
a letter after Gateaux's death.\footnote{Postcard dated from 7 December
1914.}

%s4 #&#
\section{In the Storm}

A serious danger of war had in fact been revealed only very late in
July 1914 in public opinion, and the French mostly received the
mobilization announcement on August 2nd with stupor. Like the majority,
Gateaux has been caught napping by the beginning of the war. He was
mobilized in the reserve as lieutenant of the 269th Infantry regiment,
member of the 70th infantry division.\vspace*{2pt} The diaries of the units engaged
in the war\footnote{They were put on-line by the French Ministry of
Defense \url{http://www.memoiredeshommes.sga.defense.gouv.fr}.} permit us to
follow Gateaux's part in the campaign in a very precise way. He was
appointed on August 6th as the head of the 2nd machine-gun section of
the 6th Brigade when the unit was formed in Domgermain, a suburb of the
city of Toul.\footnote{Gateaux used headed notepaper from the \textit{Hotel
\& Caf\'e de l'Europe} in Toul for his last letter to Volterra on
August 25th.} The regiment paused beyond Nancy the next day and was
supposed to go further East, but the German army's fire power stopped
it brutally a few days later near Buissoncourt, 15 kilometers east of
Nancy. At the end of August, the main task of the 70th infantry
division was to defend Nancy's southeast sector.

The centennial year 2014 was an occasion for many people to better
realize how horrific the first few weeks of the war were on the French
side. August 1914 was the worst month of the whole war in terms of
casualties, and some of the figures defy belief. On 22 August 1914, for
example, the most bloody day of the whole war for the French, 27{,}000
were killed in the French ranks (\citett{Becker2004}). The appallingly
high number of casualties was due to an alliance between the
vulnerability of the French uniform [with the famous \textit{garance}
(red) trousers up to 1915$\ldots$], the self-confidence of the
headquarters who had little consideration for their men's lives, and
the clear inadequacy of many leaders in the field. Prochasson\footnote{\citet{Prochasson2004}, pages 672--673.} advances two hypotheses to
explain why the casualties among the Grandes \'Ecoles' students (\'Ecole
Normale Sup\'erieure in particular) were so dramatic. As they were
often subordinate officers, the young students were the first killed,
as their rank placed them in the front of their section. But also, they
were sometimes moved by a kind of stronger patriotic feeling that may
have driven them to a heroism beyond their simple duty.\footnote{Prochasson mentions the famous example of Charles P\'eguy and the less
well-known one of the anthropologist Robert Hertz who unceasingly asked
his superiors for a more exposed position and was killed in April
1915.} This is evident in Marbo's testimony about her adopted son
Fernand, who explained to her that, as a socialist involved in the
fight for the understanding between peoples and peace, he wanted to
\textit{be sent on the first line in order to prove that he was as brave
as anyone else},\footnote{\^Etre envoy\'e en premi\`ere ligne afin de
prouver qu' (il \'etait) aussi courageux que n'importe qui.} and added
that {\it those who would survive will have the right to speak loudly
in front of the shirkers}.\footnote{Ceux qui survivront auront le droit
de parler haut devant les embusqu\'es (\citett{Marbo1967}, page 166).}

Gateaux's last letter to Volterra is dated August 25th. Gateaux alluded
there to the ambiguous situation of Italy. Though officially allied to
the Central Empires, the country had carefully proclaimed its
neutrality, an interesting point described at length in \citet{Rusconi2005}. Senator Volterra immediately sided with France and Great
Britain and wrote passionate letters to his French colleagues as early
as the beginning of August to express the hope that Italy would join
them.\footnote{See \citet{MazliakTazzioli2009} where Volterra's attitude
is thoroughly studied.} On 24 October 1914, in a letter to Borel, he
asked for news
\begin{quote}
from Mr. Gateaux, Mr. P\'er\`es, Mr. Boutroux and Mr.
Paul L\'evy and other young French friends. I have received a letter
from Mr. Gateaux from the battlefield and then no other. And this is
why I am very worried about his fate and that of the others.\footnote{M. Gateaux, M. P\'er\`es, M. Boutroux, M. Paul L\'evy et d'autres jeunes
amis fran\c cais. ($\ldots$) J'avais re\c cu une lettre de M. Gateaux du
champ de bataille et ensuite je n'en ai re\c cu pas d'autre c'est
pourquoi je suis tr\`es inquiet sur son compte ainsi que sur les autres.}
\end{quote}

Borel answered Volterra's letter on November 4, telling him that P\'er\`es and Boutroux were discharged and that he did not know where Gateaux
was.\footnote{The tone of this letter was slightly less confident than
the previous ones. This was the moment when the enormous losses of the
first weeks began to filter through. Borel wrote that at the \'Ecole
Normale, several young men with a bright scientific future had already
disappeared and that the responsibility of \textit{those who wanted this
war} was really terrible.} As we have seen, Gateaux was in Lorraine at
the end of August. The French army went steadily backward, and was
closer and closer to being crushed between the two wings of the German
army (one coming from the north through Belgium, the other from the
east through Lorraine and Champagne). Then occurred the unexpected \textit{miracle} of the Battle of the Marne (6--13 September 1914), which
suddenly stopped the German advance, rendering the Schlieffen Plan a
failure. Vitry-le-Fran\c cois had been occupied by the Germans during
the night of the 5th of September, but they were compelled to leave and
to withdraw toward the East on September 11th.\footnote{A vivid account
of this moment was written after the war by a witness (\citett{Nebout1922}). Though Gonthiez and Janet wrote in \citet{AnnuaireENS}
that they could \textit{easily imagine all the pain he \textup{(}Gateaux\textup{)} would
have felt when he learned that the enemy had taken the city of
Vitry-le-Fran\c cois where his poor mother had remained}, it is not
clear whether Gateaux had learnt the fact at all, due to the general
confusion. I refer to \citet{Becker2004} or to several articles of
\citet{encyclopedie2004} for the description of this phase of the war.} From
September 13th, the French went again slowly toward the East, chasing
after the retreating Germans.

At the end of September, the French and British and the German
headquarters became aware of the impossibility of any further decisive
motion on the front line running from the Aisne to Switzerland; each
realized that the only hope was to bypass their enemy in the zone
between the Aisne and the sea which was still free of soldiers.

General Joffre decided to withdraw from the Eastern part of the front
(precisely where Gateaux was) a large number of divisions and to send
them \textit{by railway} to places in Picardie, then in Artois and finally
to Flanders to try to outrun the Germans. The so-called \textit{race for
the sea} lasted two months and was very bloody.

The 70th division was transported between September 28th and October
2nd from Nancy to Lens, a distance of almost 500 km.\footnote{According to the diary of the 269th Infantry regiment, the order to board the
trains, received on September 28th, was carried out the next day. With
an impressive organizational efficiency, the trains followed a
circuitous route to join Artois: Troyes, Versailles, Rouen before
stopping at Saint-Pol sur Ternoise on October 1.} Gateaux's division
received the order to defend the East of Arras. On October 3rd,
Gateaux's regiment was in Rouvroy, a small village, 10 km southeast
from Lens, and Gateaux was killed at one o'clock in the morning while
trying to prevent the Germans from entering the village. In the
confusion of the bloodshed, the corpses were not identified before
being collected and hastily buried in improvised cemeteries. Gateaux's
body was buried near St. Anne Chapel in Rouvroy, a simple cross without
inscription marking the place.\footnote{According to the army file,
Ren\'e's mother was informed on October 4 that her son was reported missing.
On March 16th 1916, her other son and only remaining child, Ren\'e's
brother Georges, was killed in the Mort-Homme before Verdun. Much
later, Ren\'e's mother passed away on 24 February 1941 in
Vitry-le-Fran\c cois, some months after having seen her city devastated
by the German invasion.}

Ren\'e's death was officially established only on 28 December
1915.\footnote{This was done based on evidence given by Henri-Auguste
Munier-Pugin, warrant officer, and Albert Garoche, sergeant, in the
269th Infantry regiment.} But it is only long after, on 8 December
1921, that Gateaux's corpse was exhumed and formally identified, and
finally transported to the necropolis of the military cemetery of the
Bietz-Neuville St Vaast.\footnote{Gateaux's grave is number 76 at
Bietz-Neuville. Gateaux's mother was informed of this fact on 5 January
1922.} The last document of the military dossier is a letter from the
Minister of War, dated 22 June 1923, informing the mayor of
Vitry-le-Fran\c cois that the Lieutenant Ren\'e-Eug\`ene Gateaux had
officially been declared \textit{Dead for France}.

The detailed chronology of how the academic world learned of Gateaux's
death is not entirely clear. As already mentioned, the Principal of
Bar-le-Duc Lyc\'ee wrote the postcard in December 1914, but it was
clearly an answer to a letter he had received.\footnote{This postcard
is, however, a decisive link between Hadamard and the papers left by
Gateaux. It was probably addressed to Hadamard or Borel, though I found
it by chance in the huge archive of Fr\'echet material in the Paris
Academy of Science. Another possibility is that the letter was
addressed to Fr\'echet who happened to know the Proviseur as well as
Gateaux well enough to have this exchange. If this hypothesis is true,
it may be Fr\'echet who recovered Gateaux's papers and transmitted them
to Hadamard. We shall see a point below that corroborates this version.}

Only on December 10th did Borel write to Volterra about Gateaux's death
(\citett{MazliakTazzioli2009}, page  47), mentioning his anxious hope
that of the dozens of pupils of the \'Ecole Normale considered as lost,
there will be at least one or two who will come back at the end of the
war. Volterra sadly answered some days later (\citett
{MazliakTazzioli2009}, page 48) and wrote that he was sure that Ren\'e
would have had a great future. The same day a telegram was sent to the
\'Ecole Normale by Volterra in the name of the Mathematical seminar in Rome.

As early as August 1915, Hadamard took the necessary steps to obtain
the award of one of the Paris Academy's prizes for Gateaux. In a letter
dated 5 August 1915 (and probably addressed to Picard as Perpetual
Secretary), Hadamard mentioned the following: Gateaux \textit{has left
very advanced research on functional calculus \textup{(}his thesis was composed
to a great extent\textup{,} and partly published in notes to the Academy\textup{),}
research for which M. Volterra and myself have a great regard}.\footnote{(Gateaux) laisse sur le calcul fonctionnel des recherches
fort avanc\'ees (sa th\`ese \'etait en grande partie compos\'ee, et
repr\'esent\'ee par des notes pr\'esent\'ees \`a l'Acad\'emie),
recherches auxquelles M.~Volterra, comme moi-m\^eme, attache un grand
prix.} At the meeting\vadjust{\goodbreak} of 18 December 1916, the Franc\oe ur prize was
awarded to Gateaux (\citett{Hadamard1916}, pages 791--792). It is
interesting to read in Hadamard's short report the following section:
\begin{quote}
[Gateaux] was following a much more audacious way, which
promised to be very fruitful, by extending the notion of integration to
the functional domain. Nobody could predict the development and the
range this new series of research would attain. This is what has been
interrupted by events.\footnote{(Gateaux) allait s'engager dans une voie
beaucoup plus audacieuse, et qui promettait d'\^etre des plus f\'econdes, en \'etendant au domaine fonctionnel la notion d'int\'egrale.
Nul ne peut pr\'evoir le d\'eveloppement et la port\'ee qui auraient pu
\^etre r\'eserv\'es \`a cette nouvelle s\'erie de recherches. C'est
elle qui a \'et\'e interrompue par les \'ev\'enements.}
\end{quote}

It is plausible that Hadamard had only superficially looked at
Gateaux's papers, since he himself was caught in the storm of events,
losing his two sons during the summer of 1916. Nevertheless, he did at
least notice that one major interest in the last period of Gateaux's
work was integration over the space of functionals. As we shall see,
this was precisely why he spoke to L\'evy about Gateaux.\vspace*{-3pt}

%s5 #&#
\section{The Mathematical Destiny}\label{destiny}

%s5.1 #&#
\subsection{L\'evy's Interest in Infinite-Dimensional Integration}
\begin{quote}
In January 1918, I was lying on a bed in a hospital, when
I suddenly thought again of functional analysis. In my early work, I~had never thought of extending the notion of an integral to spaces with
infinite dimensions. It suddenly appeared to me that it was possible to
attack this problem starting with the notion of mean in a sphere of the
space of square summable functions.
Such a function can be approximated by a step function, the number $n$
of its distinct values growing constantly. The desired mean may then be
defined as the limit of the mean in a sphere of the $n$-dimensional
space. Obviously, this limit may not exist; but in practice, it does
often exist  (\citett{Levy1970}, page 58).
\end{quote}

Thus, L\'evy described how he became interested in infinite-dimensional
integration. It is not easy to decide whether this happened as suddenly
as he wrote, just following the train of his thoughts. Regardless, it
is sometimes forgotten today that L\'evy, before becoming one of the
major specialists in Probability theory of the 20th Century, had been a
brilliant expert in functional analysis.\footnote{On that topic, see, in
particular, \citet{Barbut2014}, pages  44--54.} As we shall see, it is a
remarkable fact that his studies in functional analysis led him rather
naturally to probabilistic formulations of problems. At the end of
1918, the Paris Academy of Sciences, following Hadamard's proposal,
decided to call upon L\'evy for the \textit{Cours Peccot} in 1919.\footnote{The Cours Peccot was (and still is) a series of lectures in
mathematics given at the Coll\`ege de France and financed by the Peccot
Foundation. It is a way to promote innovation in research by offering
financial support and an audience to a young mathematician. Borel had
been the first lecturer in 1900, followed by Lebesgue. In L\'evy's
time, the age of the lecturer was meant to be less than thirty.
However, the losses of the war had been so heavy among young men that
the choice of the thirty-three year old L\'evy was reasonable. It is
also plausible to think that Gateaux would have been a natural Peccot
lecturer had he survived the war. As L\'evy's appointment is almost
concomitant with Hadamard asking to take care of Gateaux's papers, it
is possible that there is a connection between the two events.}
L\'evy's book \textit{Le\c cons d'Analyse Fonctionnelle} (\citett{Levy1922}),
on which I shall comment later, is based on these Peccot lectures.

The first document in which the question is explicitly mentioned is a
letter to Volterra written in the early days of 1919:
\begin{quote}
As I was recently interested in the question of the
extension of the integral to functional space, I spoke about the fact
to Mr.~Hadamard who mentioned the existence of R. Gateaux's note on the
theme. But he could not give me the exact reference and I cannot find
it. $(\ldots)$ Though I am still mobilized, I am working on lectures I
hope to give at the Coll\`ege de France on the functions of lines and
equations with functional derivatives, and on this occasion I would
like to develop several chapters of the theory. ($\ldots$) I think that
the generalization of the Dirichlet problem must present greater
difficulties. Up to now, I was not able to extend your results on
functions of the first degree and your extension of Green's formula.
This is precisely due to the fact that I do not possess a convenient
expression for the integral.\footnote{M'\'etant occup\'e r\'ecemment de
la question de l'extension de la notion d'int\'egrale multiple \`a
l'espace fonctionnel, j'en ai parl\'e \`a M.~Hadamard qui m'a signal\'e
l'existence d'une note de R.~Gateaux sur ce sujet. Mais il n'a pas pu
m'en donner la r\'ef\'erence exacte et je ne puis r\'eussir \`a la
trouver. ($\ldots$) Quoiqu'encore mobilis\'e, je travaille \`a pr\'eparer
un cours que j'esp\`ere professer au Coll\`ege de France sur les
fonctions de lignes et les \'equations aux d\'eriv\'ees fonctionnelles
et \`a cette occasion, je voudrais d\'evelopper davantage certains
chapitres de la th\'eorie. ($\ldots$) Je crois que la g\'en\'eralisation
du probl\`eme de Dirichlet doit pr\'esenter plus de difficult\'es. Je
n'ai pu jusqu'ici profiter pour le cas g\'en\'eral de vos travaux sur
les fonctions du premier degr\'e et l'extension de la formule de Green.
Ceci tient pr\'ecis\'ement \`a ce que je n'ai pas encore mis la notion
d'int\'egrale multiple sous une forme commode pour ce but. (L\'evy to
Volterra, 3 January 1919.)}
\end{quote}

As can be seen from this quotation, L\'evy's views on
infinite-dimensional integration were related to his studies in
potential theory. The central problem of the classical mathematical
potential theory is to find a harmonic function $U$ in a domain $R$
with given values on the boundary $S$ (Dirichlet problem) or given
values of the normal derivatives on $S$ (Neumann problem). In 1906,
Hadamard (\citeyear{Hadamard1906}) proposed to make use of variational
techniques from Volterra's theory of functions of lines in order to
study more general forms of these problems, for instance, when the
border is moving with time, and, in particular, to find Green functions
used in the integral representation of the solutions. These problems
would make up L\'evy's thesis, defended in 1911.

As L\'evy wrote to Volterra, to study these questions in
infinite-dimensional functional spaces, one needs to be able to
integrate over these spaces. Volterra was not the only person L\'evy
had contacted. He wrote to Fr\'echet on the same topic at the very end
of the year 1918.\footnote{See \citet{Barbut2014}, page 69.} Fr\'echet
had indeed proposed
in Fr\'{e}chet  (\citeyear{Frechet1915})
a theory of integration over abstract spaces in
1915, usually considered as the first attempt to define a general
integral.\footnote{These are Kolmogorov's terms in \citet
{Kolmogorov1933}. On this matter, see, for instance, \citet{ShaferVovk2006}.}

On 6 January 1919, L\'evy wrote to Fr\'echet
\begin{quote}
About Gateaux's papers, I learned precisely yesterday
that M. Hadamard had put them in security at the \'Ecole Normale during
the war and had just taken them back. Nothing is therefore yet
published.\footnote{\citet{Barbut2014}, Lettre 2.}
\end{quote}

From this, I infer that Fr\'echet mentioned Gateaux's papers to L\'evy,
probably because he had an idea of what they contained. This could also
be a hint that the papers arrived to Hadamard during the war via Fr\'
echet and that Fr\'echet was the addressee of the postcard from the
Principal of Bar-le-Duc.

On January 12, L\'evy sent another letter to Volterra:
\begin{quote}
M. Hadamard has just found several of Gateaux's unpublished
papers at the \'Ecole Normale. I have not seen them yet but maybe I'll
find what I am looking for in them.\footnote{M.~Hadamard vient de trouver
plusieurs m\'emoires non publi\'es de Gateaux \`a l'\'Ecole Normale. Je
ne les ai pas encore vus mais peut-\^etre y trouverais-je ce que j'y
recherche.}
\end{quote}
Volterra answered on January 15, writing that none of Gateaux's
publications concerned integration. He nevertheless added
\begin{quote} Before he left Rome, we had discussed about his general
ideas on the subject, but he did not publish anything. I suppose that
in the manuscripts he had left, one may probably find some notes
dealing with the problem. I am happy that they are not lost and that
you have them in hand. The question is very interesting.\footnote{Nous avons caus\'e avant son d\'epart de Rome des id\'ees g\'en\'erales sur
ce sujet mais il n'a rien publi\'e l\`a-dessus. Je pense que dans les
notes manuscrites qu'il a laiss\'ees, on pourra bien probablement
trouver quelques notes sur ce sujet. Je suis heureux qu'elles ne soient
pas perdues et qu'elles se trouvent dans vos mains. La question est tr\`
es int\'eressante.}
\end{quote}

As already mentioned, Hadamard entrusted L\'evy with the posthumous
edition of Gateaux's papers.
He published it in three parts as Gateaux  (\citeyear{Gateaux6}, \citeyear{Gateaux7})
and (\citeyear{Gateaux8}).
In February 1919, L\'evy began to describe
the precise content of what he had found in Gateaux's papers to Fr\'echet.

%s5.2 #&#
\subsection{Gateaux's Integration of Functionals}

Integration over infinite-dimensional spaces was certainly the most
important subject considered by Gateaux. This can be read in Hadamard's
comment that follows:
\begin{quote} The fact that he chose functional calculus reveals a
broad mind, scornful of small problems or of the easy application of
known methods. But the event proved that Gateaux was able to consider
such a study under its widest and most suggestive aspect. And it is
what he indeed did, with integration over the functional field, to
speak only about this\vadjust{\goodbreak} example, the most important, that represents a
path that is new and the theory.\footnote{Le fait qu'il ait choisi le
calcul fonctionnel r\'ev\'elait un esprit aux vues larges, d\'edaigneux
du petit probl\`eme ou de l'application facile de m\'ethodes connues.
Mais le fait prouva que Gateaux \'etait capable de consid\'erer une
telle \'etude sous son aspect le plus large et le plus suggestif. Et
c'est effectivement ce qu'il fit, avec l'int\'egration sur le champ
fonctionnel, pour ne mentionner que cet exemple, le plus important, qui
repr\'esente une voie enti\`erement nouvelles et de tr\`es grandes
perspectives pour la th\'eorie (\citett{AnnuaireENS}, page 138).}
\end{quote}
Gateaux's views on integration are the subject of the first paper
edited by L\'evy in 1919 (\citett{Gateaux6}). L\'evy completed this
presentation (and considerably extended it) in Part III of \citet{Levy1922}, Chapter II, page~274.

As said before, when I commented on Gateaux's letter to Volterra
expositing his research program, Gateaux's interest in
infinite-dimensional integration originated in an attempt to extend
Cauchy's formula and his first idea was to use a Riemann-type approach.

Gateaux considered the ball\footnote{To fit better with modern
terminology, I use the word \textit{ball}, though Gateaux and L\'evy
systematically use \textit{sphere}.} consisting of all square integrable
functions over $[0,1]$ with the property $\int_0^1x(\alpha)^2\,d\alpha
\le
R^2$.\footnote{In fact, Gateaux started from a \textit{continuous}
function $x$. However, as L\'evy explained to Fr\'echet in a long
letter dated 16 February 1919 (Letter 5 in \citett{Barbut2014}), it is
more natural to consider measurable functions, that is, to work with the
(now) usual space $L^2$. This is what he does in \citet{Levy1922}.} He
defined a function $x$ to be {\it simple of order $n$} if it assumes
constant values $x_1,x_2,\ldots,x_n$ over each subinterval $[0, \frac{1}{n}[,\ldots,[\frac{n-1}{n}, 1]$. In order that a simple function $x$
belongs\vspace*{1pt} to the ball, one must therefore have $x_1^2+x_2^2+\cdots
+x_n^2\le nR^2$. The set of simple functions of order $n$ belonging to
the ball is called the {$n$th section of the ball}. This set
corresponds to a ball in $\mathbb{R}^n$ centered at 0 with radius
$\sqrt n R$.

As the volume $V_n$ of a ball with radius $\sqrt n R$ in dimension $n$
is asymptotically equivalent\vspace*{1pt} to $ {(2\pi \mathbf{e})^{n/2}\over
\sqrt{n\pi}}R^n$ (\citett{Levy1922}, page  265), it tends to zero or
infinity for $n \to\infty$, depending on the value of $R$. This fact
constitutes the central problem for the definition of the integral: in
functional space, a subset has generally a volume equal to zero or
infinity, and this forbids the direct extension of the Riemann integral
through an approximating step-function sequence.

Gateaux seems to have been the first to propose a natural way to bypass
the problem by defining the integral as a limit of mean values.
Consider a functional $U$ defined and continuous on the ball $\int_0^1x(\alpha)^2\,d\alpha\le R^2$. Its restriction $U_n$ to the $n$th
section can be considered as a continuous function of the $n$ variables
$x_1,x_2, \dots, x_n$ and, therefore, it admits a mean value
\[
\mu_n=\frac{\int_{x_1^2+x_2^2+\cdots+x_n^2\le nR^2}U_n(x_1,\ldots,
x_n)\,dx_1\cdots \,dx_n}{V_n}.
\]
Under some circumstances, the sequence $(\mu_n)$ admits a limit which
is called the mean value of $U$ over the ball of the functional space.
Gateaux's main achievement in \citet{Gateaux6} was to obtain the value
of the mean for important types of functionals.

He began by considering functionals of the type $U\dvtx  x\mapsto f[x(\alpha
_1)]$ where $x$ is a point of the functional space, $f$ a continuous
real function and $\alpha_1$ a fixed point in $[0,1]$. As $\alpha_1$ is
fixed, $x(\alpha_1)$ is one of the coordinates when $x$ is taken in the
$n$th section.\footnote{Gateaux considers this functional although it is
clearly not continuous. Gateaux had not sorted out the role of
continuity in his work on the infinite-dimensional. It is likely that
he would have improved the apparent incoherence in a subsequent
rewriting of the paper. We shall see that L\'evy fixed the question in
\citet{Levy1922}.}

Therefore, the [($n-1$)-dimensional] volume of the intersection of the
ball of radius $R$ with the plane $x(\alpha_1)=z$ (with $0\le z^2\le
nR^2$ or, equivalently, $-\sqrt n R \le z\le\sqrt n R$) is given by
\[
\bigl(\sqrt{nR^2-z^2}\bigr)^{n-1}\cdot
V_{n-1},
\]
where $V_{k}$ is the volume of the unit ball in dimension $k$. A
classical result is that for any $k\ge2$, $V_k$ satisfies the
induction formula $V_k=2V_{k-1}\int_0^{\pi/2}\cos^k \theta
\,d\theta$.

Now, the mean of the functional $U$ over the $n$th section is given by
\begin{eqnarray*}
&& \frac{1}{(\sqrt n R)^n\cdot V_n}\\[-2pt]
&& \quad{}\cdot\int_{-\sqrt n R}^{+\sqrt n R}f(z) \bigl(\bigl(
\sqrt{nR^2-z^2}\bigr)^{n-1}\cdot V_{n-1}
\bigr)\,dz.
\end{eqnarray*}
Performing the change of variables $z=R\sqrt n \theta$ transforms the
previous expression into
\[
\frac{1}{\int_{-{\pi/2}}^{\pi/2}\cos^n \theta \,d\theta
}\int_{-{\pi/2}}^{\pi/2}f(R\sqrt n \sin
\theta)\cos^n \theta \,d\theta.
\]
It is seen that the preponderant values for $\theta$ in the last
integral are those around 0, and $ \int_{-{\pi/2}}^{\pi/2}\cos^n \theta \,d\theta$ is known to be asymptotically
equivalent to $\sqrt{\frac{2\pi}{n}}$. Under ``some regularity
conditions'' for $f$, the previous expression is therefore
approximately equal to
\[
\frac{1}{\sqrt{(2\pi)/n}}\int_{-{\alpha\sqrt n }}^{\alpha\sqrt n
}f\biggl(R\sqrt n
\sin{\psi\over\sqrt n} \biggr)\cos^n {\psi\over\sqrt n}
{d\psi\over\sqrt n}
\]
for any\vadjust{\goodbreak} $\alpha>0$ and sufficiently large $n$.

Using a Taylor expansion, and letting $n$ go to infinity, the latter
expression converges to
%
%e1 #&#
\begin{equation}
\label{limit} \frac{1}{\sqrt{2\pi}}\int_{-\infty
}^{+\infty}f(R
\psi)\mathbf{e}^{-{\psi^2/2}}\,d\psi,
\end{equation}
defined by Gateaux as the mean of $U$ over the ball of all square
integrable functions over $[0,1]$ such that $\int_0^1x(\alpha)^2\,d\alpha
\le R^2$. He asserted that this result can be generalized for
functionals of the type
\begin{eqnarray*}
U(x) &=& \int_0^1d\alpha_1\cdots
\int_0^1d\alpha_p\\
&&{}\cdot f\bigl[x(
\alpha_1 ),\ldots, x(\alpha_p),
\alpha_1,\ldots, \alpha_p\bigr]
\end{eqnarray*}
for which the mean value is given by
%
%e2 #&#
\begin{eqnarray}
&& \frac{1}{(2\pi)^{p/2}}\nonumber\\
&&\quad{}\cdot\int_0^1d
\alpha_1\cdots\int_0^1d\alpha
_p\int_{-\infty}^{+\infty}dx_1\cdots
\int_{-\infty}^{+\infty
}dx_p
\nonumber
\\[-8pt]
\label{limit2}\\[-8pt]
\nonumber
&&\qquad{}\cdot f(Rx_1,
\ldots ,Rx_p,\alpha_1 , \ldots, \alpha_p)
\\
&&\quad\qquad{}\cdot\mathbf{e}^{-{(x_1^2+\cdots+x_p^2)/2}}.\nonumber
\end{eqnarray}

The rigorous existence of the limit was not explained by Gateaux, as
L\'evy wrote to Fr\'echet in his letter of 12 February 1919. Obviously,
for Gateaux, as L\'evy himself wrote in the foreword of \citet{Gateaux6}, the present state of his papers was certainly not a final
one.\footnote{As can be seen, Gateaux used a technique close to
Laplace's method for the estimation of the limit. This method for
asymptotic estimation of integrals was currently taught to students in
Paris, but usually without much care for the convergence conditions.
This may also explain that Gateaux did not pay much attention to this
aspect of the question in his manuscript.} And in the long note L\'evy
added at the end of the article (\citett{Gateaux6}, page 67), he
described the attempts made by Gateaux to obtain the limit in several
situations. For L\'evy, the priority was to fill the gap left by
Gateaux and to try to obtain the existence of the mean value for the
most general functionals.

Gateaux (\citeyear{Gateaux6}, page 52) also considered continuous (with
respect to uniform norm) functionals $U$ satisfying the following
property: for any $\varepsilon>0$, there is an $n_0$ such that, for
$n\ge n_0$ and for any two functions $x$ and
$y$ satisfying $\int_0^1x(\alpha)^2\,d\alpha\le R^2$ and assuming the
same mean value over each subinterval $[\frac{i-1}{n}, \frac{i}{n}]$,\footnote{Which is to say that $\int_{(i-1)/n}^{i/n}
x(t)\,dt=\int_{(i-1)/n}^{i/n} y(t)\,dt$ for every $i$ such that
$1\le i\le n$.} one has $\vert U(x)-U(y) \vert< \varepsilon$. Following
Gateaux, for such a functional, the mean value is given by the value at
the center 0 of the ball (the function constantly equal to 0), and it
can therefore be considered as a \textit{harmonic} functional. The
previously mentioned property of $U$ was natural to Gateaux: he had
proved in \citet{Gateaux2} that, under such a condition, a continuous
functional $U$ can be well approximated over the ball $\int_0^1x(\alpha
)^2\,d\alpha\le R^2$ by $U(y_n)$, where $y_n$ belongs to the $n$th
section of the sphere and takes on the interval $[\frac{i-1}{n},
\frac{i}{n}]$ the value $\int_{(i-1)/n}^{i/n} x(t)\,dt$.\vspace*{1pt}

After he began to scrutinize Gateaux's paper, L\'evy became convinced
that Gateaux's requirement of continuity with respect to the uniform
norm for a functional $U$ was in fact much too restrictive. As early as
16 February 1919,\footnote{\citet{Barbut2014}, page 115.} he mentioned
the fact to Fr\'echet. And in the final version of his ideas on the
question, in \citet{Levy1922}, page 277, he arrived at a striking
conclusion: under very general assumptions, such a continuous
functional takes \textit{almost everywhere} the same constant value $b$,
meaning that for any $\varepsilon>0$, the volume of the subset of
these functions $x$ in the $n$th section of the ball satisfying $\vert
U(x)-b \vert>\varepsilon$ tends to 0 with $n\to\infty$. The mean of
such a
functional is therefore obviously equal to this value $b$. L\'evy gives
(\citett{Levy1922}, page 275) a simple example illustrating this
situation. Consider the functional defined on the ball $\int_0^1x(\alpha)^2\,d\alpha\le R^2$ by $U(x)=\varphi(r)$, where $\varphi$
is a given continuous function on $\mathbb{R}_+$ and $r^2= \int_0^1x^2(\alpha)\,d\alpha$. The volume of the ball $B_n(\sqrt n R)$ with
radius $\sqrt n R$ centered in 0 in $\mathbb{R}^n$ is proportional to
$(\sqrt n
R)^n$; hence, for any given $0<\varepsilon<1$, the quotient of the
volumes of $B_n( (1-\varepsilon)\sqrt n R) )$ and $B_n(\sqrt n R)$
tends to 0, which means that when $n$ grows, the volume is more and
more concentrated close to the surface. Therefore, $\varphi(R)$ is
essentially the only value assumed by $\varphi$ in the ball counting
for the calculation of the mean.\footnote{The concentration of measure
phenomenon became an important field of research following Milman's
systematic study of asymptotic geometry in Banach spaces during the
1970s. It has many important applications, especially in probability
theory by providing exponentional inequalities of Gaussian type. See
\citet{Ledoux2001} for a panoramic view of this question.}

%s5.3 #&#
\subsection{L\'evy's Probabilistic Interpretation}

I have already mentioned that in 1919, L\'evy had his first contact
with probability theory when he was asked to teach probability at the
\'Ecole Polytechnique.\footnote{For more details about this story, I refer
the reader to \citet{BarbutMazliak2008}.} This was exactly the same
period he was studying Gateaux's papers and preparing their
publication. One may observe that probability theory takes no part in
the various notes presented by L\'evy to the Paris Academy of Sciences
as he progressed in his work on Gateaux (\citeauthor{Levy1919a}, \citeyear{Levy1919a},
\citeyear
{Levy1919b}, \citeyear{Levy1919c}, \citeyear{Levy1921}).\footnote{L\'evy began,
however, to work on independent probabilistic questions at the same
time. See, in particular, \citet{Fischer2011}, page~218 \textit{and seq.}
for L\'evy's investigations on characteristic functions and the central
limit theorem, and \citet{Barbut2014}, pages 40--44, more specifically
about L\'evy's investigations on stable distributions.} But when he
wrote his book \citet{Levy1922} he often adopted probabilistic
reasonings as relevant for his considerations about the mean in a
functional space and it seems that a kind of extraordinary junction
occurred during these years in L\'evy's mind, resulting in unifying his
mathematical interests in functional calculus and probability
theory.\footnote{Recall here his own mention that he \textit{was prepared
by functional calculus for the study of functions with an infinite
number of variables and \textup{(}that\textup{)} many of \textup{(}his\textup{)} ideas in functional
analysis became without effort ideas which could be applied in
probability} (\citett{Barbut2014}, page 156).}

Let us try to understand how probability entered L\'evy's
considerations about the mean in functional spaces [third part of \citet{Levy1922}].
Consider (\citett{Levy1922}, page 266) a given hyperplane
$H$ containing 0 in $\mathbb{R}^n$ and define the coordinate $z$ as the
distance to $H$. Let us consider the fraction of the ball centered at 0
with radius $R\sqrt n$, comprised between the hyperplanes $z=R\xi_1$
and $z=R\xi_2$. The ratio of the volume of this fraction to the total
volume of the ball is equal to
\[
\int_{\xi_1/\sqrt n}^{\xi_2/\sqrt n}\cos^n\theta \,d\theta
\over\int_{-\pi/2}^{+\pi/ 2}\cos^n\theta \,d
\theta
\]
which tends to
%
%e3 #&#
\begin{equation}
\label{gauss} \frac{1}{\sqrt{2\pi}}\int_{\xi_1}^{\xi
_2}
\mathbf{e} ^{-{x^2/2}}\,dx.
\end{equation}

More generally, consider $p$ hyperplanes containing 0 and call
$z_1,z_2,\ldots, z_p$ the distances to these hyperplanes. The volume of
the intersection of $p$ regions $R\xi_i'<z_i<R\xi_i^{\prime\prime}$
($i=1,2,\ldots, p$) is a fraction of the total volume equal to
\begin{eqnarray*}
&& \frac{1}{(2\pi)^{p/2}}\int_{\xi_1^{\prime}}^{\xi_1^{\prime\prime}}dx_1
\int_{\xi_2^{\prime}}^{\xi
_2^{\prime\prime}}dx_2\cdots\int_{\xi_p^{\prime}}^{\xi_p^{\prime\prime}}dx_p\\
&&\quad{}\cdot\mathbf{e}^{-({x_1^2+x_2^2+\cdots+x_p^2})/{2}}.
\end{eqnarray*}
This is, writes L\'evy, a direct consequence of the independence of the
random variables $z_i$, each following a Gaussian distribution
according to the previous result. In order to prove the desired
independence, writes L\'evy, it is sufficient to prove that the
conditions $z_i=R\xi_i, i=1,2,\ldots, p-1$ do not influence the
distribution of $z_p$. The intersection of these conditions is a
hyperspace $H$ with dimension $n-p+1$, included in a hyperplane $r=kR$
($r$ being the distance between 0 and $H$). Now, the intersection of
$H$ and the ball of radius $R\sqrt n$ is a ball with dimension $n-p+1$
and radius $R\sqrt{n-k^2}$, asymptotically equivalent to $R\sqrt
{n-p+1}$ when $n$ tends to infinity. Moreover, $n-p+1$ tends to
infinity with $n$. Therefore, concludes L\'evy, the distribution of
$z_p$ is given by the formula (\ref{gauss}), hence the desired independence.

As the reader can see, L\'evy's proof is based on a kind of intuitive
approach which would become his typical trademark in numerous later
works in probability. In particular, the sketchy use of conditional
densities seems almost sloppy for a modern mathematician's eye, but
L\'evy was never embarrassed with such technicalities in his proofs. For
L\'evy, the essential task was to understand the deep nature of the
mathematical situation. In so doing, he had a lot in common with
Poincar\'e's conception of what is a rigorous proof in mathematics. Not
only beyond the purely logical proofs, mathematically insignificant,
but also beyond the analytical proofs which logically deduce theorems
from definitions and axioms, Poincar\'e defended the necessity of a
specific intuition for a mathematician, a {\it geometrical} spirit
using his senses and his imagination in order to perceive this \textit{touch of something which realizes the unity of the proof}.\footnote{\label{rigor}\textit{Ce je ne sais quoi qui fait l'unit\'e de la d\'emonstration}. On Poincar\'e's conceptions, see the recent paper
(\citett{Kebaili2014}). L\'evy's intuitive approach is the precise aspect that
explains what It\^o wrote later, about his difficult work to \textit{translate} L\'evy.
\textit{At that time}, writes It\^o, \textit{it was
commonly believed that L\'evy's works were extremely difficult, since L\'evy\textup{,} a pioneer in the new mathematical field\textup{,} explained probability
theory based on his intuition. I~attempted to describe L\'evy's ideas
using precise logic that Kolmogorov might use} (\citett{Ito1998}).}

The probabilistic framework allowed L\'evy to explain Gateaux's formula
(\ref{limit}) for the mean of the functional $U(x)=f[x(\tau)]$ in what
seems to him a more convincing way (\citett{Levy1922}, page 278). If $x$
is in the ball with radius $R\sqrt n$, the probability of the event
$R\xi_1\le x(\tau) \le R\xi_2$ tends to $ {1\over\sqrt
{2\pi}}\int_{\xi_1}^{\xi_2}\mathbf{e}^{-{\xi^2/2}}\,d\xi$
when $n \to\infty
$, so that the mean of $U$ is given by (\ref{limit}). Moreover, the
mean of $U(x)=\varphi(x(t_1),x(t_2), \dots, x(t_p))$ is immediately
obtained using the fact that $x(t_1),x(t_2), \ldots, x(t_p)$ are i.i.d.
variables having a centered Gaussian distribution with variance $R^2$
(\citett{Levy1922}, page~281). Probabilistic reasoning also enables us
to explain the concentration of the mass at the surface of a ball in
the functional space (\citett{Levy1922}, page  283). By the law of large
numbers, $ \frac{x_1^2+\cdots+x_n^2}{n}$ tends to $R^2$ and,
therefore, for any $\varepsilon>0$, the probability that
$\sqrt{ \frac{x_1^2+\cdots+x_n^2}{n}}$ does not belong to
$[R-\varepsilon, R+\varepsilon]$ tends to 0 when $n\to\infty$.
Therefore, concludes L\'evy, the part of the $n$th section one must
take into account for the computation of the mean of a functional is in
the neighborhood of the surface of the sphere with radius $R\sqrt n$.

In Chapter VI (\citett{Levy1922}, Part Three, page  421), L\'evy studies
the general question of the existence of the mean for a functional. As
we have seen in the previous subsection, L\'evy considered continuity
with respect to the uniform norm as too strong a condition because it
implied that the functional is almost surely constant. In this chapter,
he highlights that in order to obtain a convenient condition for the
existence of the mean, it is necessary to look at the probability
distribution of the values of the function $x$ rather than at the
values themselves.

As a basic example he considers the mean of the functional $U(x)=F(f)$
in the ball with radius $R$, where $f$ is the probability distribution
function (called by L\'evy \textit{fonction sommatoire}) of $x$ over the
space $[0,1]$ equipped with Lebesgue measure $\lambda$.\footnote{This is
to say that $x$ is considered as a random variable on the probability
space $[0,1]$ with Lebesgue measure $\lambda$. Hence, $f(\xi) =
\lambda\{ t\in[0,1], x(t)\le\xi\}$.} L\'evy's reasoning is as
follows. If $x$ belongs to the $n$th section of the ball, it is a
function constant in each interval $[\frac{i-1}{n}, \frac{i}{n}]$ with
value $x_i$, such that $x_1^2+x_2^2+\cdots+x_n^2\le nR^2$. In the limit
$n\to\infty$, the $x_i$ are independent Gaussian random variables with
variance $R^2$, and the probability distribution function associated
with this $x$ is the Gaussian distribution function with variance $R^2$
denoted by $\varphi$.\footnote{To explain this in modern terms, consider
a sequence of independent random variables $(X_n)_{n\ge1}$, each with
the standard normal distribution. By the law of large numbers, the
sequence $\frac{1}{n} \sum_{k=1}^n \ind_{X_k\le x}$ tends almost surely
to $P(X_1\le x)$. Choose a $\omega$ for which the convergence occurs
and, for each $n$, define a random variable $Z_n$ on the probability
space $([0,1], \lambda)$ by $Z_n(t)=X_i (\omega)$ if $\frac{i-1}{n}\le
t< \frac{i}{n}$. Then $(Z_n)_{n\ge1}$ converges in distribution to the
standard normal distribution. L\'evy is extremely elliptic in his proof
(he only mentions \textit{des raisonnements connus de calcul des
probabilit\'es}). He may have had the intuition that the dependence on
$\omega$ in the previous construction would not create real
difficulties as results from the Glivenko--Cantelli theorem.} This
allows L\'evy to conclude (\citett{Levy1922}, page 424) that the mean of
$U$ is equal to $F(\varphi)$.

As a generalization of the previous result, L\'evy studies functionals
$U$ satisfying a condition which, though weaker than continuity with
respect to uniform topology, guarantees a good approximation of the
functional by its values on the $n$th section. The most general
property [called $\mathcal{H}$ by L\'evy (\citeyear{Levy1922}), page
424]
he considers is the following: for each given $\varepsilon>0$, there
is a $n$ such that, if $x$ and $y$ are two functions in the ball such
that in every interval $[\frac{i-1}{n},\frac{i}{n}]$ the probability
distribution function of $x$ and $y$ is the same,\footnote{This means
that considered as random variables on the probability space
$[\frac{i-1}{n},\frac{i}{n}]$ with probability measure $n\cdot\lambda$,
the two functions $x$ and $y$ restricted to this interval have the same
distribution.} $\vert U(y)-U(x) \vert<\varepsilon$. However, L\'evy was
not able to prove the approximation result he was looking for in all
the desired generality, but he asserted that the result was
\textit{reliable} for the functional satisfying the property $\mathcal{H}$
(\citett{Levy1922}, page~427).

As it is seen, probability reasoning is omnipresent in the Third Part
of \citet{Levy1922}. L\'evy was certainly conscious of the profound
originality of his approach and desired to convince everyone of its
interest. The complicated relations between the prominent French
mathematicians (Borel and Hadamard in the first place) and probability
theory was considered in several studies (see \citett{Bru2003} and \citett
{DurandMazliak2011} and the references included for more details). It
was observed that from the very beginning of his interest in
probability, L\'evy felt himself unjustly despised for his
choice,\footnote{On that topic, see, in particular,
\citet{BarbutMazliak2008b}.} though he was comforted by Wiener's reaction to
his approach (I shall come back on that point in the next subsection).

This lack of interest of the leading French mathematicians in
probability (Borel was the exception) may be an explanation why
absolutely no reference to probability can be located in Gateaux's
papers, even when he observed the remarkable appearance of the Gaussian
distribution in the limit expression (\ref{limit}).
In \citet{Borel1906},\footnote{Reprinted as \textit{Note I} in his book
(\citett{Borel1914}).} Borel had proved that if $B_n$ is the ball of $\mathbb{R}^n$
centered in 0 with radius $R\sqrt n$, and $V_n (u)$ the volume of the
portion $u\le x_1\le u+du$ of $B_n$, the ratio of $V_n(u)$ to the total
volume of $B_n$ tends to $\frac{1}{\sqrt{2\pi}R} \mathbf{e}^{-{u^2/(2R^2)}}\,du$.\footnote{The result, usually
known today under the name \textit{Poincar\'e's lemma}, has in fact
nothing to do with Poincar\'e, according to Diaconis and Freedman
(\citeyear{Diaconis1987}).
Moreover, Stroock (\citeyear{Stroock1994}) discovered that
Mehler had already obtained the result in 1866 in a purely analytical
context [see \citet{Stroock1994}, page 68, footnote 3 for an exact
reference and comments].} Borel's interest was statistical mechanics,
more precisely, for Maxwell and Boltzmann's kinetic theory of gases. In
his presentation, the spheres represent surfaces in the phase space of
equal total kinetic energy. In a complement to his translation of
Ehrenfests' paper on statistical mechanics in \textit{Encyclop\'edie des
Sciences Math\'ematiques} (\citett{Borel1914b}, page 273), Borel mentions
studies about the $n$-dimensional sphere as the first example of
mathematical research inspired by statistical mechanics. He even
audaciously asserts that one should consider the results about surfaces
and volumes in high dimensions as connected to statistical mechanics.
However, in contrast to Maxwell, who, in his fundamental paper in 1860,
had emphasized the coincidence between the distribution law for the
speeds of the particles and the distribution \textit{governing the
distribution of errors among observations by use of the so-called}
least-squares method,\footnote{\citet{Maxwell1860}, Prop.~IV and
following comments.} Borel did not mention any possible connection
with the law of errors in \citet{Borel1906}. The only reference is in
\citet{Borel1914}, page  66, without any probabilistic interpretation,
just mentioning that the \mbox{Gaussian} distribution function was a
well-tabulated distribution function which allows it to be used for
computations.

It is probably the desire to explain to a large audience why
probabilistic tools were useful that prompted L\'evy to write a
nontechnical paper for the \textit{Revue de M\'etaphysique et de Morale}
(\citett{Levy1924}). L\'evy explains there the general ideas leading to
his conception of the mean value, based on probability considerations
over general sets.\footnote{Interestingly, L\'evy asserts (\citett{Levy1924}, page 149) that the article is the development of his last
lecture of the Cours Peccot of 1919, meaning that the aforementioned
junction between probability and his studies in functional calculus
appeared quite early in his mind. This is corroborated by his first
letters to Fr\'echet (\citett{Barbut2014}, Letters 1--5, before February
1919). If probability is never mentioned explicitly there, one may
observe how gradually L\'evy is closer to probabilistic reasoning. A
good example is found in Letter 3 (\citett{Barbut2014}, page~55) where
L\'evy writes about his desire to find a way of expressing that functions
$u$ such that $\int u'^2$ is large are \textit{less probable}.} As an\vadjust{\goodbreak}
elementary example, he considers the situation of non-negative integers
as today in probabilistic number theory. If $f$ is a function defined
on $\mathbb{N}$ ($f$ could typically be the indicator of a subset
$A\subset\mathbb{N}$), the mean of $f$ is defined as the limit of $\frac{1}{N}\sum_{k=1}^Nf(k)$ when $N$ tends to infinity. In particular, $P(A)=\lim_{N\rightarrow+\infty}\frac{1}{N} \operatorname{Card} \{ n\in\mathbb{N}, n\in
A\}$.\footnote{Therefore, if one randomly draws a point from $\mathbb{N}$,
there is,
for instance, one chance over two that it is an even integer, a rather
comforting result for the mind$\ldots$} The paper includes a
presentation of Gateaux's work on infinite-dimensional integration and
the idea behind the extension to more general functionals. L\'evy was
probably rather satisfied with the picture he had provided in his
paper, as he decided to reprint it as an appendix in his treatise of
probability published the next year (\citett{Levy1925}). Another attempt
to disseminate his considerations on functional analysis was also done
in 1924. Henri Villat asked L\'evy to write a small booklet for his new
series \textit{M\'emorial des Sciences Math\'ematiques}. \citet{Levy1925b}
contains 56 pages and appears in fact as a survey of the book (\citett{Levy1922}). L\'evy updated his bibliography and Daniell's and Wiener's
works were now quoted.

%s5.4 #&#
\subsection{Wiener Measure: Daniell Versus Gateaux's~Integrals}

As mentioned in the \hyperref[sec1]{Introduction}, it is well beyond the scope of this
article to provide a complete description of the fundamental works
where Wiener built the first mathematical model of Brownian motion; on
that topic, I refer the reader to It\^o's comments in \citet{WienerOC}
and to \citet{Chatterji1993}, \citet{Kahane1998} and \citet{Barbut2014},
pages 54--60.
The aim of this section is more modest: to try to explain how Wiener
became acquainted with Gateaux's approach to integration and how he
eventually used it in his epoch-making paper (\citett{Wiener1923}).

In the second half of the 1910s, the British mathematician Percy J.
Daniell (1889--1946), then holding a position at the Rice Institute in
Houston, Texas, was interested in extending Lebesgue integration\vspace*{1pt} to
infinite-dimensional spaces.\footnote{A very complete description of
Daniell's work and personality can be found in the paper (\citett{Aldrich2007}).} Daniell wrote two important papers
(\citeauthor{Daniell1917}, \citeyear{Daniell1917}, \citeyear{Daniell1918}) on the subject. His approach was to consider the
integral\vadjust{\goodbreak} as an operator on functions satisfying certain properties,
such as linearity and a monotone convergence theorem on a restricted
class of functions $T_0$, and to prove that these properties allow one
to extend integration to the class $T_1$ of limits of sequences in
$T_0$. It can be seen that such a construction is directly inspired by
Lebesgue.\footnote{L\'evy always coolly accepted nonconstructive
approaches, which, for him, probably did not sufficiently reveal the
\textit{touch of something} (in the words of Poincar\'e, see note  \ref{rigor} above) at the heart of a mathematical concept. Thus, he did not
hide his moderate appreciation of Daniell's work on integration to Fr\'echet. He wrote to him {\it if nothing important has escaped me\textup{,}
Daniell has given not a definition of the integral but an extension of
the notion of integral from a restricted domain to a larger one. That
is a Lebesgue-kind work \textup{(}\citett{Barbut2014}\textup{,} page \textup{86}\textup{)}.}}

Wiener's first work on functionals (\citett{Wiener1920}) appeared in 1920.
Wiener proved there that Daniell's method can be applied to define the
integral of a functional, taking as basis $T_0$ a set of step functions
for which the integral is defined as a mean. Probably shortly before
publication, Wiener added the following footnote (\citett{Wiener1920}, page  67):
\begin{quote}
The use of mean instead of integral is found in the
posthumous papers of Gateaux (Bulletin de la Soci\'et\'e Math\'ematique
de France, 1919). This was however unknown to me at the time I wrote
this article.
\end{quote}
We do not know exactly when Wiener was informed of the existence of
Gateaux's works. A possible hypothesis is that he became aware of them
during his journey in France in 1920 when he came to the Strasbourg
International Congress and met Fr\'echet and Volterra.

The next year, Wiener published his first papers on Brownian motion. In
the first one (\citett{Wiener1921a}), he starts from Einstein's result:
at time $t$ the probability that the position $f(t)$ of a particle on a
line belongs to the interval $[x_0,x_1]$ has the form $
\frac{1}{\sqrt{\pi ct}}\int_{x_0}^{x_1}\mathbf{e}^{-{x^2/ct}}\,dx$
where $c$ is a
constant (taken equal to 1 by Wiener, corresponding to a good choice of
units). The path $x=(f(t), 0\le t\le1)$ of the particle is a
real-valued continuous function on $[0,1]$. Thus, if we consider as
functional a function of this path, a natural question arises of
defining its average value. Due to the independence of increments in
the Brownian motion, asserts Wiener, it is reasonable to associate to a
functional of the form $F=\Phi(f(t_1), \ldots, f(t_n))$ depending only
on the values of $f$ at some finite number\vadjust{\goodbreak} of values of $t$, a~mean,
denoted $A[F]$ by Wiener, defined by
\begin{eqnarray*}
&& A[F] = \frac{1}{\sqrt{\pi^nt_1(t_2-t_1)\cdots(t_n-t_{n-1})}}\cdots
\\
&&\quad{}\cdot\int_{-\infty}^{+\infty}\int_{-\infty}^{+\infty}
\cdots\int_{-\infty}^{+\infty} \Phi(x_1,\ldots,
x_n)\\
&&\qquad{}\cdot\mathbf{e}^{-{x_1^2/t_1}-{(x_2-x_1)^2/
(t_2-t_1)}-\cdots-{(x_n-x_{n-1})^2/(t_n-t_{n-1})}}\\
&&\qquad\quad{}\cdot dx_1\cdots
dx_n.
\end{eqnarray*}
In particular, observes Wiener, if $F(f)=f(t_1)^{m_1}\cdot\break 
f(t_n)^{m_n}$, one may compute an explicit value for $A[F]$. Therefore,
if a functional $F$ is analytical in the sense of Volterra, which means
that it can be expanded as a sum of functionals of the type
\begin{eqnarray*}
&& \int_0^1\cdots\int_0^1f(x_1)
\cdots f(x_n)\\
&&\hspace*{34pt}\quad{}\cdot\varphi_n(x_1, \ldots,
x_n)\,dx_1\cdots \,dx_n,\!\!\footnotemark[80]
\end{eqnarray*}
\footnotetext[80]{Wiener
considers in fact a generalization of this situation where the functionals are
defined by means of Stieltjes integrals.}%
\noindent the mean of $F$ is defined as the sum of the corresponding terms
\[
\int_0^1\cdots\int_0^1A[F_n]
\varphi_n(x_1, \ldots, x_n)\,dx_1
\cdots \,dx_n
\]
[where $F_n$ is the functional $f\mapsto f(x_1)\cdots f(x_n)$] when this
series is convergent. Wiener's paper proves that, with this definition,
the mean satisfies the classical properties of integrals such as
linearity or the possibility of exchanging infinite summation and
integration. Wiener quotes Gateaux (\citett{Wiener1921a}, Note 1, page~260) for having proposed using analytical functionals in the definition
of the mean of a functional. As we have seen, it is true that Gateaux
had such an idea in mind from the very beginning (see his programmatic
letter to Borel), but, contrary to Wiener's assertion, the idea does
not seem to be explicit in \citet{Gateaux6}. Wiener adds that Gateaux's
definition is, however, not well adapted to the treatment of Brownian motion.

\setcounter{footnote}{80}

Wiener published his second study (\citett{Wiener1921b}) in the next
issue of the Proceedings of the National Academy of Sciences. The aim
of this new paper was to show that the use of the definition of the
mean provided in \citet{Wiener1921b} allowed one to obtain a direct
proof (moreover, under a somehow lighter hypotheses) of Einstein's
formula for the mean quadratic displacement of the Brownian particle in
a viscous medium. Once again, Gateaux is mentioned as having proposed
another construction of the mean:
\begin{quote}
To determine the average value of a functional, then seems
a reasonable
problem, provided that we have some convention as to what constitutes a
normal distribution of the functions that form its arguments. Two
essentially different discussions have been given on this matter: one,
by Gateaux, being a direct generalization of the ordinary mean in
$n$-space; the other, by the author of this paper, involving
considerations from the theory of probabilities  (\citett{Wiener1921b},
page 295).
\end{quote}

During the Summer of 1922, Wiener came again to France and met L\'evy
for the first time during his vacation in Pougues les Eaux, a spa in
central France, and discussed L\'evy's book on functional analysis. L\'evy narrates the meeting in his autobiography, where he emphasizes that
Wiener was almost the only one who immediately recognized the depth of
Part III of his book [\citeauthor{Levy1922} (\citeyear{Levy1922}, \citeyear{Levy1970}, page  86---and
also on page 65)]. He adds he had reasons to think that this third part
was the origin of Wiener's memoir (\citett{Wiener1923}) on Brownian motion.

Indeed, in the introduction of \citet{Wiener1923}, Wiener pays full
tribute to L\'evy:
\begin{quote}
The present paper owes its inception to a conversation
which the author had with Professor L\'evy in regard to the relation
which the two systems of integration in infinitely many dimensions---that of L\'evy and that of the author---bear to one another. For this
indebtedness the author wishes to give full credit (\citett{Wiener1923},
page 132).
\end{quote}

Gateaux is now clearly treated by Wiener only as a precursor, and L\'evy has become the major source of inspiration. Besides, Wiener wrote
(\citett{Wiener1923}, page 132) that Gateaux had begun investigations on
integration in infinitely many dimensions which had been \textit{carried
out by L\'evy} in \citet{Levy1922}.\hskip.2pt%
\footnote{%
It took some time for
Gateaux--L\'evy or Daniell considerations on infinite-dimensional
integration to be widely known. For instance, in 1930, the Danish
mathematician B{\o}rge Jessen (1907--1993) defended a doctoral thesis
with the title \textit{Contribution to the theory of the integration of
the functions of an infinity of variables} and was totally unaware of
the previous works on the topic. See \citet{BruEid2009}.}

In \citet{Wiener1923}, Wiener reconsidered the results of his previous
papers on Brownian motion. Contrary to what he had done in \citet{Wiener1921a} where the mean of a functional $F=\Phi(f(t_1), \ldots,
f(t_n))$ of the trajectory was given \textit{a priori}, he now used L\'evy's studies of the $n$-dimensional sphere and the Gateaux--L\'evy
definition of the mean as a limit of the means over the $n$th sections
in order to:
\begin{longlist}[(3)]
\item[(1)] \textit{deduce} that at time $t$, the probability distribution of the
position is Gaussian,\footnote{\citet{Wiener1923}, pages  136--137. This
was a decisive step forward with respect to \citet{Wiener1921a} where
Wiener took Einstein's Gaussian form as a starting point}

\item[(2)] \textit{define} the related measure on the space of continuous
functions (Wiener measure),

(3) \textit{prove} the value of the mean of the aforementioned functional
he had postulated in his previous works,\footnote{\citet{Wiener1923},
page~153.}

(4) \textit{derive} the expression of the mean of an analytic functional
with a new proof.\footnote{\citet{Wiener1923}, page~165.}
\end{longlist}
Section~10 of \citet{Wiener1923} is devoted to proving that, for the
functionals previously considered, Daniell's extension of the mean
Wiener had introduced in \citet{Wiener1920} gives the same value to the
integral.\footnote{The construction of the Wiener measure via Daniell's
extension is tightly related to the theorem of extension Kolmogorov
would provide 10 years later in his \textit{Grundbegriffe} (\citett{Kolmogorov1933}). On that topic, consult \citet{ShaferVovk2006}, in
particular, Section~5.1, page~87.}

Finally, observe that Wiener's paper is not absolutely conclusive about
the use of Daniell's versus Gateaux--L\'evy's approach, though I can
certainly interpret Wiener's choice to write the paper starting from
the latter as recognition of its more intuitive character. Besides, it
is well known that L\'evy was never a great supporter of abstract
constructions of Brownian motion. In his autobiography (\citett{Levy1970}, page~98), L\'evy, who was not shy about emphasizing his
missed opportunities, regretted how he let Wiener get ahead of him in
the construction of Brownian motion though all the necessary material
was in \citet{Levy1922}. L\'evy did sometimes slightly exaggerate his
own role [as, e.g., when he wrote about Kolmogoroff's \textit{Grundbegriffe} (\citett{Levy1970}, page~68)]. In the case of Brownian
motion, however, one can understand his regrets.

The geometric approach to Brownian motion was quite fertile in the
20th century. McKean (\citeyear{McKean1973}) has explained how thinking of
the Wiener measure as a uniform distribution over the \textit{infinite-dimensional sphere of radius $\sqrt{\infty}$}, a direct
consequence of L\'evy's considerations in \citet{Levy1922}, was
successfully used by Japanese mathematicians in the 1960s to describe
the geometry of Brownian motion. In another direction, in 1969,
Gallardo (\citeyear{Gallardo1969}) made the observation that \textit{Poincar\'e's lemma} could be connected with the fact that if
$X^n(t)=(X_1(t),\ldots, X_n(t))$ is an $n$-dimensional Brownian motion
starting at 0, if one denotes by $T_n$ the first passage time of $X^n$
on the sphere centered at 0 and with radius $\sqrt n$, then
$T_n\rightarrow1$ in probability and $X^n(T_n)$ follows the uniform
distribution on the $n$-dimensional sphere of radius $\sqrt{n}$. Yor
later developed these considerations (see \citett{Yor1997}).

%s6 #&#
\section{Conclusion}

It has often been said that after World War I, the French Grandes
\'Ecoles, the \'Ecole Normale especially, were crowded with the ghosts of
the students from the 1910s who disappeared during the conflict. Of
course, these dead of the Great War were essentially very young men who
had scarcely finished their graduate studies and whose names are hardly
known to us today. Ren\'e Gateaux, who died at the age of 25 in October
1914, is an example both representative and exceptional of the student
victims of the war---exceptional because, despite being very young, he
left scientific work that could be carried on by others.

Bourbaki, when he eventually added some words about probability theory
in the chapter devoted to integration in nonlocally compact spaces of a
late edition of his \textit{El\'ements d'histoire des math\'ematiques}
(\citett{Bourbaki1984}, pages 299--302),\footnote{The complicated story of
Bourbaki's attitude to integration and, in particular, of Dieudonn\'e's
resistance to abstract integration is well known and presented in
detail in \citet{Schwartz1997}. I shall not make further comments on it
here.} mentioned the path linking Borel's consideration on kinetic
theory of gases to the Wiener measure with Gateaux's and L\'evy's works
as fundamental steps.

Though uncompleted, Gateaux's mathematical studies were recovered and
extended by Paul L\'evy for whom they became a catalyst for a renewal
of his scientific interests in probability. It is due to L\'evy's work
of editing and extension that today we remember Gateaux.
\vskip1cm

\section*{Acknowledgments}
I want to express my deep thanks to the people
who kindly contributed to the preparation of this paper. First, to
M. Pierre Gateaux, a distant relative of Ren\'e Gateaux, for his warm
welcome in Vitry and his help in my research for Ren\'e's familial
background, and to Gilbert Maheux, a local historian in Vitry le Fran\c
cois, for his enthusiastic help.   In the course of this research, I
have visited several archives. I would like to thank Giorgio Letta of
the mathematical department in Pisa and member of the Accademia dei
Lincei for his efficient help, and also those responsible for the
archives in Rome for their help and kindness. Florence Greffe and the
staff of the archives of the Academie des Sciences in Paris helped me
in obtaining permission to consult Gateaux's manuscripts deposited at
the Acad\'emie. Valuable assistance was provided by Francoise
Dauphragne at the archives of the \'Ecole Normale Sup\'erieure from whom
I obtained a photograph of the 1907 pupils. And, last but not least,
through a visit to the archives of the Paris Archdiocese and the help
of Father Ph. Ploix, I discovered the only picture I know of where
Gateaux is explicitly identified. I thank also my colleagues Olivier
Gu\'edon, Bernard Locker, St\'ephane Menozzi and Marc Yor for several
interesting discussions relating to Gateaux's life and work. I~owe to
John Aldrich and Daniel Denis the transformation of the original
version of the text into what can be called real English. Finally, I am
glad to mention the anonymous referees and to thank them for having
generously provided numerous suggestions which, I hope, have eventually
largely improved the paper.

% imsref loaded by daiva.urboniene, 2015-02-03 14:39:49
% imsref loaded by daiva.urboniene, 2015-02-04 11:06:01
% imsref loaded by daiva.urboniene, 2015-02-04 12:55:25
% imsref loaded by daiva.urboniene, 2015-02-06 08:14:53

%\begin{appendix}
%\section{}
%\end{appendix}

%\section*{Acknowledgments}

%\begin{supplement}[id=suppA]
%\slink[doi]{10.1214/00-STSXXXXSUPP} %[doi,text={...}] - jei reikia
%suskaldyti doi
%\sfilename{stsXXXX\_supp.pdf}
%\end{supplement}

%\begin{thebibliography}{99}
%\bibitem[\protect\citeauthoryear{}{}]{r1}
%\bibitem{r1}
%\end{thebibliography}
\end{document}